\documentclass[11pt]{article}

\usepackage{amssymb}
\usepackage{amsfonts}
\usepackage{amsmath}
\usepackage{amsthm}
\usepackage[margin=1.00in]{geometry}
\usepackage{enumerate}
\usepackage{amstext}
\usepackage{layout}

\numberwithin{equation}{section}

\newtheorem{theorem}{Theorem}[section]
\newtheorem{corollary}{Corollary}[section]
\newtheorem{lemma}{Lemma}[section]

\newtheorem{remark}{Remark}[section]

\begin{document}

\noindent {\bf\large{A comparison theorem for the law of large numbers in Banach spaces}}

\vskip 0.3cm

\noindent {\bf Deli Li\footnote{Deli Li, Department of Mathematical
Sciences, Lakehead University, Thunder Bay, Ontario, Canada P7B
5E1\\ e-mail: dli@lakeheadu.ca} $\cdot$ Han-Ying Liang\footnote{Department of Mathematics, 
Tongji University, 1239 Siping Road, Shanghai 200092, China\\ hyliang83@yahoo.com}
\footnote{Corresponding author: Han-Ying Liang (Telephone: 86-21-65983242, FAX: 86-21-65983242)}}

\vskip 0.3cm

\noindent {\bf Abstract}~~Let $(\mathbf{B}, \|\cdot\|)$ be a real separable Banach space. Let
$\{X, X_{n}; n \geq 1\}$ be a sequence of i.i.d. {\bf B}-valued random variables and set $S_{n} =
\sum_{i=1}^{n}X_{i},~n \geq 1$. Let $\{a_{n}; n \geq 1\}$ and $\{b_{n}; n \geq 1\}$ be increasing 
sequences of positive real numbers such that $\lim_{n \rightarrow \infty} a_{n} = \infty$ and 
$\left\{b_{n}/a_{n};~ n \geq 1 \right\}$ is a nondecreasing sequence. 
In this paper, we provide a comparison theorem for the law of large numbers for i.i.d. 
{\bf B}-valued random variables. That is, we show that 
$\displaystyle \frac{S_{n}- n \mathbb{E}\left(XI\{\|X\| \leq b_{n} \} \right)}{b_{n}} \rightarrow 0$ almost surely 
(resp. in probability) for every {\bf B}-valued random variable $X$ with
$\sum_{n=1}^{\infty} \mathbb{P}(\|X\| > b_{n}) < \infty$ 
(resp. $\lim_{n \rightarrow \infty}n\mathbb{P}(\|X\| > b_{n}) = 0$) 
if $S_{n}/a_{n} \rightarrow 0$ almost surely (resp. in probability) 
for every symmetric {\bf B}-valued random variable $X$ with 
$\sum_{n=1}^{\infty} \mathbb{P}(\|X\| > a_{n}) < \infty$ 
(resp. $\lim_{n \rightarrow \infty}n\mathbb{P}(\|X\| > a_{n}) = 0$). To establish this 
comparison theorem for the law of large numbers, we invoke two tools: 
1) a comparison theorem for sums of independent {\bf B}-valued random variables and, 
2) a symmetrization procedure for the law of large numbers for sums of independent 
{\bf B}-valued random variables. A few consequences of our main results are provided. 

~\\


\vskip 0.3cm


\vskip 0.3cm

\noindent {\bf Running Head}: Law of large numbers

\section{Introduction and the main results}
Let $(\Omega, \mathcal{F}, \mathbb{P})$ be a probability space and
let $(\mathbf{B}, \| \cdot \| )$ be a real separable Banach space 
equipped with its Borel $\sigma$-algebra $\mathcal{B}$ 
($=$ the $\sigma$-algebra generated by the class of open subsets of
$\mathbf{B}$ determined by $\|\cdot\|$). A {\bf B}-valued random variable
$X$ is defined as a $\mathcal{B}$-measurable function from $(\Omega, \mathcal{F})$
into {\bf B}. Let $ \{X, X_{n};~n \geq 1 \}$ be a sequence of independent and identically
distributed (i.i.d.) {\bf B}-valued random variables and put $S_{n} = \sum_{i=1}^{n} X_{i}$,
$n \geq 1$. The strong law of large numbers (SLLN) and the weak law of large numbers (WLLN)
for i.i.d. {\bf B}-valued random variables have been studied by many authors. 
The classical Kolmogorov SLLN in real separable Banach spaces was
established by Mourier (1953). The extension of the Kolmogorov-Marcinkiewicz-Zygmund SLLN 
to {\bf B}-valued random variables is independently due to Azlarov and Volodin (1981) and de
Acosta (1981).

\vskip 0.3cm

\noindent {\bf Theorem A}.~{\rm (Azlarov and Volodin (1981) and de Acosta (1981))}. 
{\it Let $0 < p < 2$ and let $\{X, X_{n}; ~n \geq 1\}$ be a sequence of i.i.d. 
{\bf B}-valued random variables. Then
\[
\lim_{n \rightarrow \infty} \frac{S_{n}}{n^{1/p}} = 0~~\mbox{almost surely (a.s.) if and only if}
~~\mathbb{E}\|X\|^{p} < \infty~~\mbox{and}~~\frac{S_{n}}{n^{1/p}}
\rightarrow_{\mathbb{P}} 0.
\]
Here and below ``$\rightarrow_{\mathbb{P}}$" stands for convergence in probability. 
}

\vskip 0.3cm

Let $\{R, R_{n}; ~n \geq 1\}$ be a {\it Rademacher sequence}; 
that is, $\{R_{n};~n \geq 1\}$ is a sequence of i.i.d. random variables 
with $\mathbb{P}\left(R = 1\right) = \mathbb{P}\left(R = -1\right) = 1/2$. Let
$\mathbf{B}^{\infty} =
\mathbf{B}\times\mathbf{B}\times\mathbf{B}\times \cdots$ and define
\[
\mathcal{C}(\mathbf{B}) = \left\{(v_{1}, v_{2}, ...) \in
\mathbf{B}^{\infty}: ~\sum_{n=1}^{\infty} R_{n}v_{n}
~~\mbox{converges in probability} \right\}.
\]
Let $1 \leq p \leq 2$. Then $\mathbf{B}$ is said to be of {\it
Rademacher type $p$} if there exists a constant $0 < C < \infty$
such that
\[
\mathbb{E}\left\|\sum_{n=1}^{\infty}R_{n}v_{n} \right\|^{p} \leq C
\sum_{n=1}^{\infty}\|v_{n}\|^{p}~~\mbox{for all}~(v_{1}, v_{2}, ...)
\in \mathcal{C}(\mathbf{B}).
\]
The following remarkable theorem, which is due tode Acosta (1981), provides 
a characterization of  Rademacher type $p$ Banach spaces.

\vskip 0.3cm

\noindent {\bf Theorem B}.~{\rm (de Acosta (1981))}. 
{\it Let $1 \leq p < 2$. Then the following two
statements are equivalent:
\begin{align*}
& {\bf (i)} \quad \mbox{The Banach space $\mathbf{B}$ is of
Rademacher type $p$.}\\
& {\bf (ii)} \quad \mbox{For every sequence $\{X, X_{n}; ~n \geq 1 \}$
of i.i.d. {\bf B}-valued variables},
\end{align*}
\[
\lim_{n \rightarrow \infty} \frac{S_{n}}{n^{1/p}} = 0~~\mbox{a.s. if
and only if}~~\mathbb{E}\|X\|^{p} < \infty~~\mbox{and}~~\mathbb{E}X
= 0.
\]
}

\vskip 0.3cm

Let $0 < p \leq 2$. Then $\mathbf{B}$ is said to be of {\it stable type $p$} if
\[
\sum_{n=1}^{\infty} \Theta_{n}v_{n} ~~\mbox{converges a.s. whenever}~~
\{v_{n}: ~n \geq 1\} \subseteq \mathbf{B} ~~\mbox{with}~~
\sum_{n=1}^{\infty} \|v_{n}\|^{p} < \infty,
\]
where $\{\Theta_{n}; ~n \geq 1 \}$ is a sequence of i.i.d. stable random variables 
each with characteristic function $\psi(t) = \exp \left\{-|t|^{p}\right\}, ~- \infty < t <
\infty$. Equivalent characterizations of a Banach space being of stable type $p$,
properties of stable type $p$ Banach spaces, as well as various
relationships between the conditions ``Rademacher type $p$" and
``stable type $p$" may be found in Maurey and Pisier (1976),
Woyczy\'{n}ski (1978), Marcus and Woyczy\'{n}ski (1979), Rosi\'{n}ski
(1980), Pisier (1986), and Ledoux and Talagrand (1991). It is well known that 
every real separable Banach space $\mathbf{B}$ is of stable type $p$ for 
all $p \in (0, 1)$ and that if {\bf B} is of stable type $p$ for some $p \in [1, 2]$, 
then {\bf B} is of Rademacher type $p$.

A remarkable characterization of stable type $p$ Banach spaces
was provided by  Marcus and Woyczy\'{n}ski (1979). Specifically,  Marcus and Woyczy\'{n}ski (1979) 
proved the following theorem.

\vskip 0.3cm

\noindent {\bf Theorem C}.~{\rm (Marcus and Woyczy\'{n}ski (1979))}. 
{\it Let $1 \leq p < 2$. Then the following two
statements are equivalent:
\begin{align*}
& {\bf (i)} \quad \mbox{The Banach space $\mathbf{B}$ is of
stable type $p$.}\\
& {\bf (ii)} \quad \mbox{For every symmetric sequence $\{X, X_{n}; ~n \geq 1 \}$
of i.i.d. {\bf B}-valued variables},
\end{align*}
\[
\frac{S_{n}}{n^{1/p}} \rightarrow_{\mathbb{P}} 0~~\mbox{if
and only if}~~\lim_{n \rightarrow \infty} n \mathbb{P}\left(\|X\| > n^{1/p}\right) = 0.
\]
}

\vskip 0.3cm

It is well known that if $\mathbf{B}$ is of stable type $p$ for some $p \in [1,
2]$, then $\mathbf{B}$ is of stable type $q$ for all $q \in (0, p)$. In other words,
by Theorem C, we have the following conclusion: 

\vskip 0.2cm

\noindent {\bf Conclusion A}.~{\it Let $1 \leq p < 2$ and $0 < q \leq p$. If for every symmetric sequence 
$\{X, X_{n}; ~n \geq 1 \}$ of i.i.d. {\bf B}-valued random variables
\[
\frac{S_{n}}{n^{1/p}} \rightarrow_{\mathbb{P}} 0 ~~\mbox{if and only if}~~\lim_{n \rightarrow \infty}
n \mathbb{P}\left(\|X\| > n^{1/p} \right) = 0,
\]
then for every symmetric sequence $\{X, X_{n}; ~n \geq 1 \}$ of i.i.d. {\bf B}-valued random variables
\[
\frac{S_{n}}{n^{1/q}} \rightarrow_{\mathbb{P}} 0 ~~\mbox{if and only if}~~\lim_{n \rightarrow \infty}
n \mathbb{P}\left(\|X\| > n^{1/q}\right) = 0.
\]
}

Motivated by Conclusion A and Theorems B and C, in this paper we will establish what we call 
a comparison theorem for the law of large numbers for i.i.d. {\bf B}-valued random variables. 
The main results of this paper are the following Theorems 1.1 and 1.2. 

\vskip 0.3cm

\begin{theorem}
{\rm{(Comparison theorem for WLLN).}}
Let $(\mathbf{B}, \|\cdot\|)$ be a real separable Banach space.
Let $\{a_{n}; n \geq 1\}$ and $\{b_{n}; n \geq 1\}$ be increasing
sequences of positive real numbers such that 
\begin{equation}
\lim_{n \rightarrow \infty} a_{n} = \infty~ \mbox{and}~ \left\{b_{n}/a_{n}; ~n \geq 1 \right\}~
\mbox{is a nondecreasing sequence}.
\end{equation}
Suppose that, for every symmetric sequence $\{X, X_{n}; ~n \geq 1 \}$ of i.i.d. {\bf B}-valued random variables,
\begin{equation}
\frac{S_{n}}{a_{n}} \rightarrow_{\mathbb{P}} 0 ~~\mbox{if and only if}~~\lim_{n \rightarrow \infty}
n \mathbb{P}(\|X\| > a_{n}) = 0.
\end{equation}
Then, for every sequence $\{X, X_{n}; ~n \geq 1 \}$ of i.i.d. {\bf B}-valued random variables, we have that
\begin{equation}
\frac{S_{n}- \gamma_{n}}{b_{n}} 
\rightarrow_{\mathbb{P}} 0 ~~\mbox{or}~~\limsup_{n \rightarrow \infty} 
\mathbb{P} \left(\frac{\left\|S_{n}- \gamma_{n} \right\|}{b_{n}} 
> \lambda \right) > 0 ~~\forall~\lambda > 0
\end{equation}
according as
\begin{equation}
\lim_{n \rightarrow \infty}
n \mathbb{P}(\|X\| > b_{n}) = 0~~\mbox{or}~~\limsup_{n \rightarrow \infty}
n \mathbb{P}(\|X\| > b_{n}) >0.
\end{equation}
Here and below $\gamma_{n} = n \mathbb{E}\left(XI\{\|X\| \leq b_{n} \} \right)$, $n \geq 1$.
\end{theorem}

\vskip 0.3cm

\begin{theorem}
{\rm{(Comparison theorem for SLLN).}}
Let $(\mathbf{B}, \|\cdot\|)$ be a real separable Banach space.
Let $\{a_{n}; n \geq 1\}$ and $\{b_{n}; n \geq 1\}$ be increasing
sequences of positive real numbers with (1.1). Suppose that, for every symmetric sequence
$\{X, X_{n}; ~n \geq 1 \}$ of i.i.d. {\bf B}-valued random variables
\begin{equation}
\lim_{n \rightarrow \infty} \frac{S_{n}}{a_{n}} = 0 
~~\mbox{a.s. if and only if}~~\sum_{n=1}^{\infty}\mathbb{P}(\|X\| > a_{n}) < \infty.
\end{equation}
Then, for every sequence $\{X, X_{n}; ~n \geq 1 \}$ of i.i.d. {\bf B}-valued random variables, we have that
\begin{equation}
\lim_{n \rightarrow \infty} \frac{S_{n}- \gamma_{n}}{b_{n}} 
= 0~~\mbox{a.s. or}~~\limsup_{n \rightarrow \infty} \frac{\left\|S_{n}- \gamma_{n} \right\|}{b_{n}} 
= \infty~~\mbox{a.s.}
\end{equation}
according as
\begin{equation}
\sum_{n=1}^{\infty}\mathbb{P}(\|X\| > b_{n}) < \infty ~~\mbox{or}~~ = \infty.
\end{equation}
\end{theorem}

\vskip 0.2cm

\begin{remark}
Under the assumptions of Theorem 1.1, we conclude from Theorem 1.1 that, for every 
sequence $\{X, X_{n}; ~n \geq 1 \}$ of i.i.d. {\bf B}-valued random variables,
\[
\frac{S_{n} - \gamma_{n}}{b_{n}} \rightarrow_{\mathbb{P}} 0 ~~\mbox{if and only if}~~
\lim_{n \rightarrow \infty}
n \mathbb{P}(\|X\| > b_{n}) = 0,
\]
\[
\limsup_{n \rightarrow \infty} 
\mathbb{P} \left(\frac{\left\|S_{n}- \gamma_{n} \right\|}{b_{n}} 
> \lambda \right) > 0 ~~\forall~\lambda > 0 ~~\mbox{if and only if}~~
\limsup_{n \rightarrow \infty}
n \mathbb{P}(\|X\| > b_{n}) > 0.
\]
Hence that
\[
\frac{S_{n} - \gamma_{n}}{b_{n}} \nrightarrow_{\mathbb{P}} x ~~\forall ~x \in \mathbf{B}\backslash \{0\}.
\]
\end{remark}

\vskip 0.2cm

\begin{remark}
Under the assumptions of Theorem 1.2, it follows from the conclusion of Theorem 1.2 that, for every 
sequence $\{X, X_{n}; ~n \geq 1 \}$ of i.i.d. {\bf B}-valued random variables,
\[
\lim_{n \rightarrow \infty} \frac{S_{n}- \gamma_{n}}{b_{n}} 
= 0~~\mbox{a.s. if and only if}~~\sum_{n=1}^{\infty}\mathbb{P}(\|X\| > b_{n}) < \infty,
\]
\[
\limsup_{n \rightarrow \infty} \frac{S_{n}- \gamma_{n}}{b_{n}} 
= \infty~~\mbox{a.s. if and only if}~~\sum_{n=1}^{\infty}\mathbb{P}(\|X\| > b_{n}) = \infty.
\]
Hence that there does not exist a sequence $\{X, X_{n}; ~n \geq 1 \}$ of i.i.d. 
{\bf B}-valued random variables such that
\[
0 < \limsup_{n \rightarrow \infty}
\frac{\left\|S_{n} - \gamma_{n} \right\|}{b_{n}} < \infty~~\mbox{a.s.}
\]
\end{remark}

\vskip 0.2cm

Combing Theorem 1.1 and Theorem C above, we immediately obtain the following two results.

\vskip 0.2cm

\begin{corollary}
Let $1 \leq p < 2$ and let $\{a_{n}; n \geq 1\}$ be an increasing sequences of positive 
real numbers such that 
\begin{equation}
\lim_{n \rightarrow \infty} a_{n} = \infty~ \mbox{and}~ \left\{n^{1/p}/a_{n}; ~n \geq 1 \right\}~
\mbox{is a nondecreasing sequence}.
\end{equation}
Let $(\mathbf{B}, \|\cdot\|)$ be a real separable Banach space such that, for every symmetric sequence 
$\{X, X_{n}; ~n \geq 1 \}$ of i.i.d. {\bf B}-valued random variables, (1.2) holds. 
Then the Banach space {\bf B} is of stable type $p$. 
\end{corollary}

\vskip 0.2cm

\begin{corollary}
Let $(\mathbf{B}, \|\cdot\|)$ be a real separable Banach space.
Let $1 \leq p < 2$ and let $\{b_{n}; n \geq 1\}$ be an increasing sequences of positive 
real numbers such that 
\begin{equation}
\left\{b_{n}/n^{1/p}; ~n \geq 1 \right\}~\mbox{is a nondecreasing sequence}.
\end{equation}
If $\mathbf{B}$ is of stable type $p$, then, for every sequence $\{X, X_{n}; ~n \geq 1 \}$ 
of i.i.d. {\bf B}-valued random variables, (1.3) and (1.4) are equivalent.
\end{corollary}

\vskip 0.2cm

From Corollary 1.1, we see that the Conclusion A is true in general.

Similarly, by Theorem 1.2 and Theorem B above, we obtain the following two results.

\vskip 0.2cm

\begin{corollary}
Let $1 \leq p < 2$. Let $\{a_{n}; n \geq 1\}$ be an increasing sequences of positive real numbers with (1.8).
Let $(\mathbf{B}, \|\cdot\|)$ be a real separable Banach space such that, for every symmetric sequence 
$\{X, X_{n}; n \geq 1 \}$ of i.i.d. {\bf B}-valued random variables, (1.5) holds. 
Then the Banach space {\bf B} is of Rademacher type p. 
\end{corollary}

\vskip 0.2cm

\begin{corollary}
Let $1 \leq p < 2$ and let $\{b_{n}; n \geq 1\}$ be an increasing sequences of positive 
real numbers with (1.9). If $\mathbf{B}$ is of Rademacher type $p$, then, for every sequence 
$\{X, X_{n}; ~n \geq 1 \}$ of i.i.d. {\bf B}-valued random variables, (1.6) and (1.7) are equivalent.
\end{corollary}

\vskip 0.2cm

If $\mathbf{B} = (-\infty, \infty)$, Theorems 1 and 2 of Feller (1946) follow
from Corollary 1.4. 

\vskip 0.2cm

The proofs of Theorems 1.1 and 1.2 will be provided in Section 4. To establish Theorems 1.1 and
1.2, we invoke two tools: 1) a comparison theorem for sums of independent {\bf B}-valued random 
variables, which is obtained in Section 2 and, 2) a symmetrization procedure for the law of 
large numbers for sums of independent {\bf B}-valued random variables, which is provided in
Section 3. These two tools are of independent interest.

\section{A comparison theorem for sums of independent {\bf B}-valued random variables}

Throughout this section, let $\{R, R_{n}; ~n \geq 1\}$ be a Rademacher sequence.
To establish Theorems 2.1 and 2.2, we invoke the remarkable contraction principle 
discovered by Kahane (1968). Some further extensions have been obtained by 
Hoffmann-J{\o}rgensen (1973, 1974, 1976). The main result of the contraction principle
is expressed in Theorem 4.4 of Ledoux and Talagrand (1991). 

The following result is the second part of Theorem 4.4 of Ledoux and Talagrand (1991).  

\begin{lemma}
Let $\{x_{n};~ n \geq 1\}$ be a {\bf B}-valued sequence and $\{\alpha_{n}; ~n \geq 1 \}$ a real-valued
sequence such that $\sup_{n \geq 1}|\alpha_{n}| \leq 1$. Then we have, for every $n \geq 1$ and all $t \geq 0$,
\[
\mathbb{P} \left(\left\|\sum_{i=1}^{n} \alpha_{i} R_{i} x_{i} \right\| > t \right) 
\leq 2 \mathbb{P} \left(\left\|\sum_{i=1}^{n} R_{i} x_{i} \right\| > t \right).
\]
\end{lemma}

\vskip 0.3cm

By using Lemma 2.1 above, we establish in Theorem 2.1 a comparison theorem for sums 
of independent {\bf B}-valued random variables.

\begin{theorem}
Let $\varphi(\cdot)$ and $\psi(\cdot)$ be two continuous and increasing functions
defined on $[0, \infty)$ such that $\varphi(0) = \psi(0) = 0$ and
\begin{equation}
\lim_{t \rightarrow \infty} \varphi(t) = \infty 
~~\mbox{and}~~\frac{\psi(\cdot)}{\varphi(\cdot)} ~\mbox{is nondecreasing function on}~[0, \infty).
\end{equation}
Here we define $\frac{\varphi(0)}{\psi(0)} = \lim_{t \rightarrow 0^{+}} \frac{\varphi(t)}{\psi(t)}$.
For $n \geq 1$, set $a_{n} = \varphi(n)$ and $b_{n} = \psi(n)$. Then we have:

{\bf (i)}~~Let $\{x_{n};~ n \geq 1\}$ be a {\bf B}-valued sequence such that $\|x_{n}\| \leq b_{n}, ~n \geq 1$. 
Then we have, for every $n \geq 1$ and all $t \geq 0$,
\begin{equation}
\mathbb{P}\left(\left\|\sum_{i=1}^{n} R_{i}x_{i} \right\| > t b_{n} \right) 
\leq 2 \mathbb{P} \left(\left\|\sum_{i=1}^{n} R_{i}
\varphi\left(\psi^{-1}(\|x_{i}\|)\right) \frac{x_{i}}{\|x_{i}\|} \right\| > t a_{n} \right).
\end{equation}
Here $\displaystyle \varphi\left(\psi^{-1}(\|0\|)\right) \frac{0}{\|0\|} \stackrel{\Delta}{=} 0$ since 
$\displaystyle \lim_{x \rightarrow 0} \varphi\left(\psi^{-1}(\|x\|)\right) \frac{x}{\|x\|} = 0$.

{\bf (ii)}~~If $\{V_{n};~n \geq 1 \}$ is a sequence of independent 
and symmetric {\bf B}-valued random variables, then we have, for every 
$n \geq 1$ and all $t \geq 0$,
\begin{equation}
\mathbb{P}\left(\left\|\sum_{i=1}^{n} V_{i} \right\| > t b_{n} \right) 
\leq 4 \mathbb{P} \left(\left\|\sum_{i=1}^{n} \varphi\left(\psi^{-1}(\|V_{i}\|)\right) 
\frac{V_{i}}{\|V_{i}\|} \right\| > t a_{n} \right) 
+ \sum_{i=1}^{n}\mathbb{P}\left(\|V_{i}\| > b_{n} \right).
\end{equation}
\end{theorem}

{\it Proof}~~Clearly we have, for every 
$n \geq 1$ and all $t \geq 0$,
\[
\begin{array}{lll}
\mbox{$\displaystyle 
\mathbb{P}\left(\left\|\sum_{i=1}^{n} R_{i}x_{i} \right\| > t b_{n} \right)$}
&=& \mbox{$\displaystyle 
\mathbb{P}\left(\left\|\sum_{i=1}^{n} \left(\frac{a_{n}}{b_{n}} 
\cdot \frac{\psi\left(\psi^{-1}(\|x_{i}\|) \right)}{\varphi\left(\psi^{-1}(\|x_{i}\|) \right)}
\right) R_{i} \frac{\varphi\left(\psi^{-1}(\|x_{i}\|)\right)}{\psi\left(\psi^{-1}(\|x_{i}\|) \right)} x_{i}
\right\| > t a_{n} \right)$}\\
&&\\
&=& \mbox{$\displaystyle \mathbb{P}\left(\left\|\sum_{i=1}^{n} \left(\frac{\varphi(n)}{\psi(n)} 
\cdot \frac{\psi\left(\psi^{-1}(\|x_{i}\|) \right)}{\varphi\left(\psi^{-1}(\|x_{i}\|) \right)}
\right) R_{i} \varphi\left(\psi^{-1}(\|x_{i}\|)\right) \frac{x_{i}}{\|x_{i}\|}
\right\| > t a_{n} \right)$}.
\end{array}
\]
Since $\varphi(\cdot)$ and $\psi(\cdot)$ are two continuous and increasing functions defined on
$[0, \infty)$ satisfying $\varphi(0) = \psi(0) = 0$ and (2.1), we see that $\varphi^{-1}(\cdot)$ 
is also a continuous and increasing function defined on
$[0, \infty)$ such that $\psi^{-1}(0) = 0$, $\lim_{t \rightarrow \infty} \psi^{-1}(t) = \infty$, and
\[
0 \leq \frac{\psi\left(\psi^{-1}(t)\right)}{\varphi\left(\psi^{-1}(t)\right)} 
\leq \frac{\psi(n)}{\varphi(n)} = \frac{b_{n}}{a_{n}}~~\mbox{whenever}~~
0 \leq t \leq \psi(n) = b_{n}.
\]
Note that $\|x_{n}\| \leq b_{n} = \psi(n)$, $n \geq 1$.
We thus conclude that, for every $n \geq 1$,
\[
0 \leq \frac{\psi\left(\psi^{-1}(\|x_{i}\|) \right)}{\varphi\left(\psi^{-1}(\|x_{i}\|) \right)}
\leq \frac{\psi(n)}{\varphi(n)} ~~\mbox{for}~ i = 1, 2, ..., n
\]
and hence that, for every $n \geq 1$, 
\[
0 \leq \frac{\varphi(n)}{\psi(n)} 
\cdot \frac{\psi\left(\psi^{-1}(\|x_{i}\|) \right)}{\varphi\left(\psi^{-1}(\|x_{i}\|) \right)}
\leq 1~~\mbox{for}~ i = 1, 2, ..., n.
\]
By applying Lemma 2.1 we thus have, for every $n \geq 1$ and all $t \geq 0$,
\[
\begin{array}{ll}
& \mbox{$\displaystyle 
\mathbb{P}\left(\left\|\sum_{i=1}^{n} \left(\frac{\varphi(n)}{\psi(n)} 
\cdot \frac{\psi\left(\psi^{-1}(\|x_{i}\|) \right)}{\varphi\left(\psi^{-1}(\|x_{i}\|) \right)}
\right) R_{i} \varphi\left(\psi^{-1}(\|x_{i}\|)\right) \frac{x_{i}}{\|x_{i}\|}
\right\| > t a_{n} \right)$}\\
&\\
& \mbox{$\displaystyle 
\leq 2 \mathbb{P} \left(\left\|\sum_{i=1}^{n} R_{i}
\varphi\left(\psi^{-1}(\|x_{i}\|)\right) \frac{x_{i}}{\|x_{i}\|} \right\| > t a_{n} \right)$}
\end{array}
\]
proving Part {\bf (i)}.

We now turn to the proof of Part {\bf (ii)}. For every $n \geq 1$, write
\[
V_{n,i} = V_{i}I\left\{\|V_{i} \| \leq b_{n} \right\},
~~T_{i} = \varphi\left(\psi^{-1}(\|V_{i}\|)\right)\frac{V_{i}}{\|V_{i}\|}, 
~~T_{n,i} = \varphi\left(\psi^{-1}(\|V_{n,i}\|)\right)\frac{V_{n,i}}{\|V_{n,i}\|}, ~~i = 1, ..., n.
\]
Clearly we have, for all $t \geq 0$,
\begin{equation}
\mathbb{P}\left(\left\|\sum_{i=1}^{n} V_{i} \right\| > t b_{n} \right) 
\leq \mathbb{P} \left( \left\|\sum_{i=1}^{n} V_{n,i} \right\| > t b_{n} \right) 
+ \sum_{i=1}^{n}\mathbb{P}\left(\|V_{i}\| > b_{n} \right).
\end{equation}
Note that
\[
\left\{\|V_{i} \| \leq b_{n} \right\} = 
\left\{\psi^{-1}(\|V_{i} \|) \leq n \right\} =
\left\{\varphi\left(\psi^{-1}(\|V_{i} \|) \right) \leq a_{n} \right\} 
= \left\{\|T_{i} \| \leq a_{n} \right\}, ~~i = 1, ..., n.
\]
Thus it is easy to see that
\[
T_{n,i} = T_{i} I\left\{\|T_{i} \| \leq a_{n} \right\}, ~~i = 1, ..., n.
\]
Since $\{V_{i};~i = 1, ..., n \}$ is a sequence of independent 
and symmetric {\bf B}-valued random variables,  $\{V_{n,i};~i = 1, ..., n \}$, $\{T_{i};~i = 1, ..., n \}$, and
$\{T_{n,i};~i = 1, ..., n \}$ are sequences of independent and symmetric {\bf B}-valued random variables.
Let $\{R, R_{n}; ~n \geq 1\}$ be a Rademacher sequence independent of $\{V_{n};~n \geq 1 \}$.
Then $\left\{R_{i}V_{n,i}; ~i = 1, ..., n \right\}$ has the same distribution as 
$\left\{V_{n,i}; ~i = 1, ..., n \right\}$ in ${\bf B}^{n}$ and 
$\left\{R_{i}T_{n,i};~i = 1, ..., n \right\}$ has the same distribution as 
$\left\{T_{n, i};~ i = 1, ..., n \right\}$ in ${\bf B}^{n}$. Since $\|V_{n,i}\| \leq b_{n}, ~i = 1, ..., n$, 
by applying (2.2), we have, for all $t \geq 0$,
\begin{equation}
\begin{array}{lll}
\mbox{$\displaystyle
\mathbb{P} \left( \left\|\sum_{i=1}^{n} V_{n,i} \right\| > t b_{n} \right)$}
&=& \mbox{$\displaystyle
\mathbb{P} \left( \left\|\sum_{i=1}^{n} R_{i} V_{n,i} \right\| > t b_{n} \right)$}\\
&&\\
&=& 
\mbox{$\displaystyle
\mathbb{E}\left(\mathbb{P}\left(\left. \left\|\sum_{i=1}^{n} R_{i} V_{n,i} \right\| > t b_{n} 
\right| V_{1}, ..., V_{n} \right) \right)$}\\
&&\\
&\leq&
\mbox{$\displaystyle 
2 \mathbb{E}\left(\mathbb{P}\left(\left. \left\|\sum_{i=1}^{n} R_{i} T_{n,i} \right\| > t a_{n} 
\right| V_{1}, ..., V_{n} \right) \right)$}\\
&&\\
&=&
\mbox{$\displaystyle 
2 \mathbb{P} \left( \left\|\sum_{i=1}^{n} R_{i} T_{n,i} \right\| > t a_{n} \right)$}\\
&&\\
&=&
\mbox{$\displaystyle 
2 \mathbb{P} \left( \left\|\sum_{i=1}^{n} T_{n,i} \right\| > t a_{n} \right)$}.
\end{array}
\end{equation}
Since $\{T_{i};~ i = 1, ..., n \}$ is a sequence of independent and symmetric {\bf B}-valued random variables,
it follows that $\left\{T_{i}I\left\{\|T_{i}\| \leq a_{n} \right\}
- T_{i}I\left\{\|T_{i}\| > a_{n} \right\}; ~i = 1, ..., n \right\}$ has the same distribution as 
$\{T_{i};~i = 1, ..., n \}$ in ${\bf B}^{n}$.
Note that 
\[
\sum_{i=1}^{n} T_{n,i} = \frac{\sum_{i=1}^{n}T_{i} + \sum_{i=1}^{n} \left(T_{i}I\left\{\|T_{i}\| \leq a_{n} \right\}
- T_{i}I\left\{\|T_{i}\| > a_{n} \right\} \right)}{2}.
\]
We thus have, for every $n \geq 1$ and all $t \geq 0$,
\begin{equation}
\begin{array}{lll}
\mbox{$\displaystyle
\mathbb{P} \left( \left\|\sum_{i=1}^{n} T_{n,i} \right\| > t a_{n} \right)$}
&\leq& \mbox{$\displaystyle
\mathbb{P} \left( \left\|\sum_{i=1}^{n} T_{i} \right\| > t a_{n} \right)$}\\
&&\\
&& \mbox{$\displaystyle
+ \mathbb{P} \left( \left\|\sum_{i=1}^{n} \left(T_{i}I\left\{\|T_{i}\| \leq a_{n} \right\}
- T_{i}I\left\{\|T_{i}\| > a_{n} \right\} \right) \right\| > t a_{n} \right)$}\\
&&\\
&=& \mbox{$\displaystyle
2 \mathbb{P} \left( \left\|\sum_{i=1}^{n} T_{i} \right\| > t a_{n} \right).$}
\end{array}
\end{equation}
Now we can see that (2.3) follows from (2.4), (2.5), and (2.6). ~$\Box$

\section{A symmetrization procedure for the law of large numbers}

Symmetrization procedure is one of the most basic and powerful tools in probability theory,
particularly in the study of the limit theorems for sums of random variables; see, for example, 
Lemma 7.1 of Ledoux and Talagrand (1991) and Li (1988). In this section a symmetrization procedure 
for the LLN for sums of independent {\bf B}-valued random variables is established in the following 
Theorem.

\begin{theorem}
Let $\{Y_{n};~n \geq 1 \}$ be a sequence of independent 
and {\bf B}-valued random variables. Let $\{Y_{n}^{\prime};~n \geq 1 \}$ be an
independent copy of $\{Y_{n};~n \geq 1 \}$. Write $\hat{Y}_{n} = Y_{n} - Y_{n}^{\prime}, ~n \geq 1$.
Let $\{a_{n}; n \geq 1\}$ be increasing 
sequences of positive real numbers such that 
\begin{equation}
\lim_{n \rightarrow \infty} a_{n} = \infty ~~\mbox{and}~~\limsup_{n \rightarrow \infty}
\mathbb{P}\left(\|Y_{n}\| > \frac{a_{n}}{2} \right) < \frac{1}{2}.
\end{equation}
Then we have the following two statements.

{\bf (i)} ~~Symmetrization procedure for WLLN: 
\begin{equation}
\frac{\sum_{i=1}^{n} Y_{i} - \sum_{i=1}^{n} \mathbb{E} \left(Y_{i} 
I\{\|Y_{i}\| \leq a_{n} \} \right)}{a_{n}} \rightarrow_{\mathbb{P}} 0
\end{equation}
if and only if
\begin{equation}
\frac{\sum_{i=1}^{n} \hat{Y}_{i}}{a_{n}} \rightarrow_{\mathbb{P}} 0.
\end{equation}

{\bf (ii)} ~~Symmetrization procedure for SLLN: 
\begin{equation}
\frac{\sum_{i=1}^{n} Y_{i} - \sum_{i=1}^{n} \mathbb{E} \left(Y_{i} 
I\{\|Y_{i}\| \leq a_{n} \} \right)}{a_{n}} \rightarrow 0~~\mbox{a.s.}
\end{equation}
if and only if
\begin{equation}
\frac{\sum_{i=1}^{n} \hat{Y}_{i}}{a_{n}} \rightarrow 0 ~~\mbox{a.s.}
\end{equation}
\end{theorem}

\vskip 0.2cm

Theorem 3.1 tells us that for studying the LLN for $\sum_{i=1}^{n} Y_{i}/a_{n},~n \geq 1$, it
is enough to study the LLN for $\sum_{i=1}^{n} \hat{Y}_{i}/a_{n},~n \geq 1$, reducing ourselves
to symmetric random variables. To our knowledge, Theorem 3.1 (especially Theorem 3.1 (i)) is new. 
The following Corollaries 3.1 and 3.2 are two consequences of theorem 3.1.

\vskip 0.2cm

\begin{corollary}
Under the assumptions of Theorem 3.1, we have 
\[
\frac{\sum_{i=1}^{n} Y_{i}}{a_{n}} \rightarrow 0 ~~\mbox{a.s. (resp. in probability)}
\]
if and only if 
\[
\lim_{n \rightarrow \infty} \frac{\sum_{i=1}^{n} \mathbb{E} \left(Y_{i} 
I\{\|Y_{i}\| \leq a_{n} \} \right)}{a_{n}} = 0 ~~\mbox{and}~~
\frac{\sum_{i=1}^{n} \hat{Y}_{i}}{a_{n}} \rightarrow 0 ~~\mbox{a.s. (resp. in probability)}.
\]
\end{corollary}

\vskip 0.2cm

\begin{corollary}
Let $\{X, X_{n}; ~n \geq 1\}$ be a sequence of i.i.d. 
{\bf B}-valued random variables. Let $\{X_{n}^{\prime};~n \geq 1 \}$ be an
independent copy of $\{X_{n};~n \geq 1 \}$. Write $S_{n} = \sum_{i=1}^{n} X_{i}$,
$S_{n}^{\prime} = \sum_{i=1}^{n} X_{i}^{\prime}$, $n \geq 1$. Let $\{a_{n}; n \geq 1\}$ 
be increasing sequences of positive real numbers such that 
$\lim_{n \rightarrow \infty} a_{n} = \infty$. Then we have 
\[
\frac{S_{n}- n \mathbb{E}\left(XI\{\|X\| \leq a_{n} \} \right)}{a_{n}}  \rightarrow 0 ~~\mbox{a.s. (resp. in probability)}
\]
if and only if 
\[
\frac{S_{n} - S_{n}^{\prime}}{a_{n}} \rightarrow 0 ~~\mbox{a.s. (resp. in probability)}.
\]
\end{corollary}

\vskip 0.2cm

\noindent 
{\it Proof of Theorem 3.1}~~We first establish Theorem 3.1(i). Since $\{Y_{n}^{\prime};~n \geq 1 \}$ is an
independent copy of $\{Y_{n};~n \geq 1 \}$, we see that for $n \geq 1$,
\[
\frac{\sum_{i=1}^{n} \hat{Y}_{i}}{a_{n}} = \frac{\sum_{i=1}^{n} Y_{i} - \sum_{i=1}^{n} \mathbb{E} \left(Y_{i} 
I\{\|Y_{i}\| \leq a_{n} \} \right)}{a_{n}} - \frac{\sum_{i=1}^{n} Y_{i}^{\prime} 
- \sum_{i=1}^{n} \mathbb{E} \left(Y_{i}^{\prime} I\{\|Y_{i}^{\prime}\| \leq a_{n} \} \right)}{a_{n}}
\]
so that (3.3) follows from (3.2). 

We now prove that (3.3) implies (3.2). Since $\{\hat{Y}_{n};~n \geq 1 \}$ is sequence of 
independent and symmetric {\bf B}-valued variables, by the remarkable L\'{e}vy inequality in a Banach space setting 
(see, e.g., see Proposition 2.3 of Ledoux and Talagrand (1991)), we have that for every $n \geq 1$,
\begin{equation}
\mathbb{P} \left(\frac{\max_{1 \leq i \leq n}\|\hat{Y}_{i}\|}{a_{n}} > t \right) 
\leq 2 
\mathbb{P} \left(\frac{\left \| \sum_{i=1}^{n} \hat{Y}_{i} \right \|}{a_{n}} > t \right)~~\forall~t \geq 0.
\end{equation}
It thus follows from (3.3) and (3.6) that 
\begin{equation}
\frac{\max_{1 \leq i \leq n}\|\hat{Y}_{i}\|}{a_{n}} \rightarrow_{\mathbb{P}} 0.
\end{equation}
Also, by independence, (3.7) implies that
\[
\mathbb{P} \left(\frac{\max_{1 \leq i \leq n}\|\hat{Y}_{i}\|}{a_{n}} > t \right) =
1 - \prod_{i=1}^{n} \left(1 - \mathbb{P} \left(\frac{\|\hat{Y}_{i}\|}{a_{n}} > t \right) \right)
\rightarrow 0 ~~\mbox{as}~~n \rightarrow \infty ~~\forall~t \geq 0
\]
which is equivalent to
\begin{equation}
\sum_{i=1}^{n} \mathbb{P} \left(\frac{\|\hat{Y}_{i}\|}{a_{n}} > t \right) \rightarrow 0 ~~\mbox{as}~~n \rightarrow 
\infty ~~ \forall ~t > 0.
\end{equation}
It easily follows from (3.1) that there exists $n_{0} \geq 1$ such that
\[
\max_{1 \leq i \leq n} \mathbb{P} \left(\|Y_{i}\| > \frac{a_{n}}{2} \right) < \frac{1}{2} ~~\forall ~n \geq n_{0}.
\]
Note that $\left \{\|Y_{i}\| \leq a_{n}/2, \|Y_{i}^{\prime}\| > a_{n} \right \} 
\subseteq \left \{ \|Y_{i} - Y_{i}^{\prime} \| > a_{n}/2 \right \} = \left \{ \|\hat{Y}_{i} \| > a_{n}/2 \right \}$,
$1 \leq i \leq n$. We thus have that 
\[
\sum_{i=1}^{n} \mathbb{P} \left(\|Y_{i}\| > a_{n} \right) 
= \sum_{i=1}^{n} \mathbb{P} \left(\|Y_{i}^{\prime}\| > a_{n} \right) 
\leq 2 \mathbb{P} \left(\|\hat{Y}_{i} \| > a_{n}/2 \right)
~~\forall ~n \geq n_{0}
\]
and hence by (3.8), it follows that
\begin{equation}
\sum_{i=1}^{n} \mathbb{P} \left(\|{Y}_{i}\| > a_{n} \right) \rightarrow 0 ~~\mbox{as}~~n \rightarrow 
\infty.
\end{equation}
Write 
\[
Y_{n,i} = Y_{i} I\{\|Y_{i}\| \leq a_{n} \}, ~Y_{n,i}^{\prime} = Y_{i}^{\prime} I\{\|Y_{i}^{\prime}\| \leq a_{n} \},
~ \hat{Y}_{n,i} = Y_{n,i} - Y_{n,i}^{\prime}, ~~i = 1, ..., n, ~n \geq 1.
\]
It is easy to see that for every $n \geq 1$,
\[
\mathbb{P} \left(\frac{\left \| \sum_{i=1}^{n} \hat{Y}_{n,i} - \sum_{i=1}^{n} \hat{Y}_{i} \right\|}{a_{n}} \neq 0 \right)
\leq 2 \sum_{i=1}^{n} \mathbb{P} \left(\|Y_{i}\| > a_{n} \right) \rightarrow 0 ~~\mbox{as}~ n \rightarrow \infty
\]
and hence by (3.3), we conclude that
\begin{equation}
\frac{\sum_{i=1}^{n} \hat{Y}_{n,i}}{a_{n}} \rightarrow_{\mathbb{P}} 0.
\end{equation}
Again by the remarkable L\'{e}vy inequality in a Banach space setting, (3.10) ensures that
\begin{equation}
\frac{\max_{1 \leq i \leq n} \|\hat{Y}_{n,i}\|}{a_{n}} \rightarrow_{\mathbb{P}} 0.
\end{equation}
Note that
\[
\frac{\max_{1 \leq i \leq n} \|\hat{Y}_{n,i}\|}{a_{n}} \leq 2 ~~\forall ~n \geq 1.
\]
It thus follows from (3.11) and Lebesgue's dominated convergence theorem that
\begin{equation}
\lim_{n \rightarrow \infty} 
\mathbb{E} \left( \frac{\max_{1 \leq i \leq n} \|\hat{Y}_{n,i}\|}{a_{n}} \right) = 0.
\end{equation}
By Proposition 6.8 of Ledoux and Talagrand (1991, p. 156), we get that for every $n \geq 1$,
\[
\mathbb{E} \left(\frac{\left\|\sum_{i=1}^{n} \hat{Y}_{n,i} \right\|}{a_{n}} \right) 
\leq 6 \mathbb{E} \left( \frac{\max_{1 \leq i \leq n} \|\hat{Y}_{n,i}\|}{a_{n}} \right)
+ 6t_{n},
\]
where
\[
t_{n} = \inf \left\{t > 0;~ \mathbb{P} 
\left( \frac{\left\|\sum_{i=1}^{n} \hat{Y}_{n,i} \right\|}{a_{n}} > t \right) \leq \frac{1}{24} \right\}.
\]
It is easy to see that (3.10) implies that $\lim_{n \rightarrow \infty} t_{n} = 0$. It thus follows from
(3.12) that
\begin{equation}
\lim_{n \rightarrow \infty} \mathbb{E} \left(\frac{\left\|\sum_{i=1}^{n} \hat{Y}_{n,i} \right\|}{a_{n}} \right) 
= 0.
\end{equation}
Since, for every $n \geq 1$, $\sum_{i=1}^{n} Y_{n,i}$ and $\sum_{i=1}^{n} Y_{n,i}^{\prime}$ are i.i.d.
{\bf B}-valued random variables, it follows from (2.5) of Ledoux and Talagrand (1991, p. 46) that 
\[
\mathbb{E} \left(\frac{\left\|\sum_{i=1}^{n} Y_{n,i} 
- \sum_{i=1}^{n} \mathbb{E} \left(Y_{i} I\{\|Y_{i}\| \leq a_{n} \} \right) \right\|}{a_{n}} \right)
\leq \mathbb{E} \left(\frac{\left\|\sum_{i=1}^{n} \hat{Y}_{n,i} \right\|}{a_{n}} \right).
\]
By (3.13), we get that
\[
\lim_{n \rightarrow \infty} \mathbb{E} \left(\frac{\left\|\sum_{i=1}^{n} Y_{n,i} 
- \sum_{i=1}^{n} \mathbb{E} \left(Y_{i} I\{\|Y_{i}\| \leq a_{n} \} \right) \right\|}{a_{n}} \right) = 0
\]
and hence that
\begin{equation}
\frac{\sum_{i=1}^{n} Y_{n,i} 
- \sum_{i=1}^{n} \mathbb{E} \left(Y_{i} I\{\|Y_{i}\| \leq a_{n} \} \right)}{a_{n}} 
\rightarrow_{\mathbb{P}} 0.
\end{equation}
Note that, for every $n \geq 1$,
\[
\begin{array}{ll}
& \mbox{$\displaystyle \frac{\sum_{i=1}^{n} Y_{i} - \sum_{i=1}^{n} \mathbb{E} \left(Y_{i} 
I\{\|Y_{i}\| \leq a_{n} \} \right)}{a_{n}}$}\\
&\\
& \mbox{$\displaystyle = \frac{\sum_{i=1}^{n} Y_{n,i} - \sum_{i=1}^{n} \mathbb{E} \left(Y_{i} 
I\{\|Y_{i}\| \leq a_{n} \} \right)}{a_{n}} 
+ \frac{\sum_{i=1}^{n} Y_{i}I\{\|Y_{i}\| \leq a_{n} \}}{a_{n}}.$}
\end{array}
\]
Thus (3.2) follows from (3.14) and (3.9).

We now establish Theorem 3.1(ii). Clearly, we only need to show that (3.4) follows from (3.5). 
Since (3.5) implies (3.3), we see that (3.2) holds. Note that $\{Y_{n}^{\prime};~n \geq 1 \}$ is an
independent copy of $\{Y_{n};~n \geq 1 \}$. It thus follows from (3.5) that
\begin{equation}
\begin{array}{ll}
& \mbox{$\displaystyle \frac{\sum_{i=1}^{n} Y_{i} - \sum_{i=1}^{n} \mathbb{E} \left(Y_{i} 
I\{\|Y_{i}\| \leq a_{n} \} \right)}{a_{n}} - \frac{\sum_{i=1}^{n} Y_{i}^{\prime} 
- \sum_{i=1}^{n} \mathbb{E} \left(Y_{i}^{\prime}
I\{\|Y_{i}^{\prime} \| \leq a_{n} \} \right)}{a_{n}} $}\\
& \\ 
& \mbox{$\displaystyle = \frac{\sum_{i=1}^{n} \hat{Y}_{i}}{a_{n}}  \rightarrow 0
~~\mbox{as}~~ n \rightarrow \infty.$}
\end{array}
\end{equation}
By applying either Theorem 3 (3.3) of Li (1988) or Lemma 7.1 of Ledoux and Talagrand 
(1991, p. 179), (3.4) follows from (3.15) and (3.2). This completes the proof of Theorem 3.1.
~$\Box$

\vskip 0.2cm

\noindent 
{\it Proof of Corollary 3.1}~~Note that, for every $n \geq 1$,
\[
\frac{\sum_{i=1}^{n} \hat{Y}_{i}}{a_{n}} = \frac{\sum_{i=1}^{n} Y_{i}}{a_{n}} 
- \frac{\sum_{i=1}^{n} Y_{i}^{\prime}}{a_{n}},
\]
\[
\frac{\sum_{i=1}^{n} \mathbb{E} \left(Y_{i} 
I\{\|Y_{i}\| \leq a_{n} \} \right)}{a_{n}} = \frac{\sum_{i=1}^{n} Y_{i}}{a_{n}} 
- \frac{\sum_{i=1}^{n} Y_{i} - \sum_{i=1}^{n} \mathbb{E} \left(Y_{i} 
I\{\|Y_{i}\| \leq a_{n} \} \right)}{a_{n}},
\]
and
\[
\frac{\sum_{i=1}^{n} Y_{i}}{a_{n}} 
= \frac{\sum_{i=1}^{n} Y_{i} - \sum_{i=1}^{n} \mathbb{E} \left(Y_{i} 
I\{\|Y_{i}\| \leq a_{n} \} \right)}{a_{n}} + \frac{\sum_{i=1}^{n} \mathbb{E} \left(Y_{i} 
I\{\|Y_{i}\| \leq a_{n} \} \right)}{a_{n}}.
\]
Thus, by Theorem 3.1, Corollary 3.1 follows immediately. ~$\Box$

\vskip 0.2cm

\noindent 
{\it Proof of Corollary 3.2}~~Since $\lim_{n \rightarrow \infty} a_{n} = \infty$
and $\{X, X_{n}; ~n \geq 1\}$ is a sequence of i.i.d. {\bf B}-valued random variables,
we have that
\[
\mathbb{P}\left(\|X_{n}\| > \frac{a_{n}}{2} \right) = \mathbb{P}\left(\|X\| > \frac{a_{n}}{2} \right)
\rightarrow 0 ~~\mbox{as}~~n \rightarrow \infty
\]
and hence that (3.1) holds for $Y_{n} = X_{n}$, $n \geq 1$. Thus Corollary 3.2 follows from Theorem 3.1 
immediately. ~$\Box$

\section{Proofs of Theorems 1.1 and 1.2}

Throughout this section, $\{a_{n}; n \geq 1\}$ and $\{b_{n}; n \geq 1\}$ are increasing
sequences of positive real numbers with (1.1). Write
\[
I(1) = \left \{i;~ b_{i} \leq 2 \right \} ~~\mbox{and}~~
I(m) = \left \{i;~2^{m-1} < b_{i} \leq 2^{m} \right\},~~m \geq 2.
\]
It follows from (1.1) that $0 < b_{n} \uparrow \infty$ and
\[
\{1, 2, 3, ..., n, ... \} = \bigcup_{m=1}^{\infty} I(m).
\]
Note that $I(m), ~m \geq 1$ are mutually exclusive sets. 
Thus there exist positive integers $k_{n}, m_{n}, ~n \geq 1$ such that
\[
k_{1} < k_{2} < ... < k_{n} < ..., ~~m_{1} < m_{2} < ... < m_{n} < ...,
\]
\[
\{1, 2, 3, ..., n, ... \} = \bigcup_{n=1}^{\infty} I\left(m_{n} \right), ~~\mbox{and}~~I\left(m_{1} \right) 
= \left\{1, ..., k_{1} \right\}, ~I\left(m_{n} \right) = \left\{k_{n-1} + 1, ..., k_{n} \right\}, ~n \geq 2.
\]

To prove Theorems 1.1 and 1.2, we use the following two preliminary lemmas. 

\vskip 0.2cm

\begin{lemma}
Then there exist two continuous and increasing functions $\varphi(\cdot)$ and $\psi(\cdot)$
defined on $[0, \infty)$ such that (2.1) holds and 
\begin{equation}
\varphi(0) = \psi(0) = 0, ~\varphi(n) = a_{n}, ~\psi(n) = b_{n}, ~ n \geq 1.
\end{equation}
\end{lemma}

\noindent {\it Proof}~~Let $a_{0} = b_{0} = 0$. Let
\[
\varphi(t) = a_{n-1} + \left(a_{n} - a_{n-1} \right) (t - n + 1), ~n - 1 \leq t < n, ~n \geq 1
\]
and
\[
\psi(t) = b_{n-1} + \left(b_{n} - b_{n-1} \right) (t - n + 1), ~n - 1 \leq t < n, ~n \geq 1.
\]
Clearly, $\varphi(\cdot)$ and $\psi(\cdot)$ are two continuous and increasing functions 
defined on $[0, \infty)$ such that (4.1) holds. We now verify that (2.1) holds with
the chosen $\varphi(\cdot)$ and $\psi(\cdot)$. Note that (1.1) implies that, 
for $n-1 < t < n$ and $n \geq 1$,
\[
\begin{array}{lll}
\mbox{$\displaystyle 
\left(\frac{\psi(t)}{\varphi(t)} \right)^{\prime} $} 
&=& 
\mbox{$\displaystyle
\frac{\psi^{\prime}(t) \varphi(t) - \psi(t) \varphi^{\prime}(t)}{\varphi^{2}(t)}$}\\
&&\\
&=& 
\mbox{$\displaystyle
\frac{a_{n-1}b_{n} - b_{n-1}a_{n}}{\varphi^{2}(t)}$}\\
&&\\
&=& 
\mbox{$\displaystyle
\frac{a_{n-1}a_{n} \left(\frac{b_{n}}{a_{n}} - \frac{b_{n-1}}{a_{n-1}} \right)}{\varphi^{2}(t)}$}\\
&&\\
& \geq &
\mbox{$\displaystyle
0, $}
\end{array} 
\]
where $b_{0}/a_{0} \stackrel{\Delta}{=} b_{1}/a_{1}$. Thus (2.1) follows. ~$\Box$

\vskip 0.2cm

\begin{lemma}
Let $\{V_{n};~n \geq 1 \}$ be a sequence of independent and symmetric {\bf B}-valued random variables.
Set $k_{0} = 0$. We have the following two statements.

{\bf (i)} ~~If 
\begin{equation}
\frac{\sum_{i=1}^{n}V_{i}}{a_{n}} \rightarrow 0~~\mbox{a.s.,} 
\end{equation}
then
\begin{equation}
\sum_{n=1}^{\infty} 
\mathbb{P}\left(\left\|\sum_{i = k_{n-1} + 1}^{k_{n}} V_{i} \right\| > \epsilon a_{k_{n}} \right)
< \infty ~~\forall~ \epsilon > 0.
\end{equation}

{\bf (ii)}
\begin{equation}
\frac{\sum_{i=1}^{n}V_{i}}{b_{n}} \rightarrow 0~~\mbox{a.s.}
\end{equation}
if and only if
\begin{equation}
\sum_{n=1}^{\infty} 
\mathbb{P}\left(\left\|\sum_{i=k_{n-1}+1}^{k_{n}} V_{i} \right\| > \epsilon b_{k_{n}} \right)
< \infty ~~\forall~ \epsilon > 0.
\end{equation}
\end{lemma}

\noindent {\it Proof}~~We first prove Part (i). Clearly, (4.2) implies that
\begin{equation}
\frac{\left\|\sum_{i = k_{n-1} + 1}^{k_{n}} V_{i} \right\|}{a_{k_{n}}} 
\leq \frac{\left\|\sum_{i = 1}^{k_{n}} V_{i} \right\|}{a_{k_{n}}} + \left(\frac{a_{k_{n-1}}}{a_{k_{n}}} \right)
\frac{\left\|\sum_{i = 1}^{k_{n-1}} V_{i} \right\|}{a_{k_{n-1}}} \rightarrow 0 ~~\mbox{a.s. as}~~n \rightarrow \infty.
\end{equation}
Since $\sum_{i=k_{n-1} + 1}^{k_{n}} V_{i}, ~n \geq 1$ are independent {\bf B}-valued random variables, (4.3) follows 
from (4.6) and the Borel-Cantelli lemma.

We now establish Part (ii). From the proof of Part (i), we only need to show that (4.4) follows from (4.5).
Since $\{V_{n};~n \geq 1 \}$ is a sequence of independent and symmetric {\bf B}-valued random variables,
by the remarable L\'{e}vy inequality in a Banach space setting (see, e.g., see Proposition 2.3 of Ledoux and Talagrand 
(1991)), we have that for every $n \geq 1$,
\[
\mathbb{P} \left( \max_{k_{n-1} < k \leq k_{n}}
\left \|\sum_{i= k_{n-1} + 1}^{k} V_{i} \right \| > \epsilon b_{k_{n}}  \right) 
\leq 2 \mathbb{P}\left(\left\|\sum_{i = k_{n-1} + 1}^{k_{n}}V_{i} \right\| > \epsilon b_{k_{n}} \right)
~~\forall~\epsilon \geq 0.
\]
Thus it follows from (4.5) that
\[
\sum_{n=1}^{\infty} \mathbb{P} \left( \max_{k_{n-1} < k \leq k_{n}}
\left \|\sum_{i= k_{n-1} + 1}^{k} V_{i} \right \| > \epsilon b_{k_{n}}  \right) < \infty 
~~\forall~\epsilon \geq 0
\]
which ensures that
\begin{equation}
A_{n} \stackrel{\Delta}{=}
\frac{\max_{k_{n-1} < k \leq k_{n}} \left \|\sum_{i= k_{n-1} + 1}^{k} V_{i} \right \|}{b_{k_{n}}}
\rightarrow 0 ~~\mbox{a.s.}
\end{equation}
Now by the Toeplitz lemma, we conclude from (4.7) that
\[
\begin{array}{lll}
\mbox{$\displaystyle 
\max_{k_{n-1} < k \leq k_{n}} \frac{\left \|\sum_{i= 1}^{k} V_{i} \right \|}{b_{k}} $}
&\leq& 
\mbox{$\displaystyle 
2 \max_{k_{n-1} < k \leq k_{n}} \frac{\left \|\sum_{i= 1}^{k} V_{i} \right \|}{b_{k_{n}}} $}\\
&&\\
&\leq&
\mbox{$\displaystyle 
2 \left(\sum_{j=1}^{n-1} \frac{\left \|\sum_{i= k_{j-1} + 1}^{k_{j}} V_{i} \right \|}{b_{k_{n}}}
+ A_{n} \right)$}\\
&&\\
&\leq&
\mbox{$\displaystyle 
2 \sum_{j=1}^{n} \left(\frac{b_{k_{j}}}{b_{k_{n}}} \right) A_{j} $}\\
&&\\
&\leq& 
\mbox{$\displaystyle 
4 \sum_{j=1}^{n} \left(\frac{2^{m_{j}}}{2^{m_{n}}} \right) A_{j} $} \\
&&\\
&\rightarrow& 
\mbox{$\displaystyle 
0 ~~\mbox{a.s.~as}~ n \rightarrow \infty,$}
\end{array}
\]
i.e., (4.4) holds. ~$\Box$

\vskip 0.2cm

With the preliminaries accounted for, Theorems 1.1 and 1.2 may be proved.

\vskip 0.2cm

\noindent {\it Proof of Theorem 1.1}~~To establish the conclusion of Theorem 1.1,
it suffices to show that, for every sequence $\{X, X_{n}; ~n \geq 1 \}$ 
of i.i.d. {\bf B}-valued random variables, the following three statements are equivalent:
\begin{equation}
\frac{S_{n}- \gamma_{n}}{b_{n}} 
\rightarrow_{\mathbb{P}} 0,
\end{equation}
\begin{equation}
\lim_{n \rightarrow \infty} \mathbb{P}
\left( \frac{\left\|S_{n}- \gamma_{n} \right\|}{b_{n}} > \lambda \right) 
= 0 ~~\mbox{for some constant}~\lambda \in (0, \infty),
\end{equation}
\begin{equation}
\lim_{n \rightarrow \infty}
n \mathbb{P}(\|X\| > b_{n}) = 0.
\end{equation}
Here and below $\gamma_{n} = n \mathbb{E}\left(XI\{\|X\| \leq b_{n} \} \right)$, $n \geq 1$.

Since (4.9) obviously follows from (4.8), it suffices to establish the implications ``(4.8) $\Rightarrow$ (4.10)", 
``(4.10) $\Rightarrow$ (4.8)", and ``(4.9) $\Rightarrow$ (4.8)". 

Since $\{X, X_{n}; ~n \geq 1 \}$ is a sequence of i.i.d. {\bf B}-valued random 
variables, arguing as in the proof of implication ``(3.2) $\Rightarrow$ (3.9)", we see that
(4.8) implies (4.10).

We now show that (4.10) implies (4.8). To see this, let $\{X, X_{n}; ~n \geq 1 \}$ 
be a sequence of i.i.d. {\bf B}-valued random variables with (4.10). Set
\[
\tilde{X} = \frac{X - X^{\prime}}{2}, ~~\tilde{X}_{n} = \frac{X_{n} - X^{\prime}_{n}}{2}, ~~n \geq 1
\]
where $\{X^{\prime}, X_{n}^{\prime};~n \geq 1 \}$ is an independent copy of $\{X, X_{n};~n \geq 1 \}$.
Clearly,
\[
\mathbb{P} (\|\tilde{X} \| > t) \leq \mathbb{P} (\|X\| > t) + \mathbb{P} (\|X^{\prime}\| > t) 
= 2 \mathbb{P} (\|X\| > t) ~~\forall ~t > 0.
\]
Thus $\{\tilde{X}, \tilde{X}_{n}; ~n \geq 1 \}$ is a sequence of i.i.d. 
symmetric {\bf B}-valued random variables such that
\begin{equation}
\lim_{n \rightarrow \infty} n \mathbb{P} \left(\|\tilde{X}\| > b_{n} \right) = 0.
\end{equation}
Since $\{a_{n}; n \geq 1\}$ and $\{b_{n}; n \geq 1\}$ are increasing
sequences of positive real numbers with (1.1), by Lemma 4.1, there exist
two continuous and increasing functions $\varphi(\cdot)$ and $\psi(\cdot)$
defined on $[0, \infty)$ such that both (2.1) and (4.1) hold. Write
\[
Y = \varphi\left(\psi^{-1}(\|\tilde{X}\|)\right) 
\frac{\tilde{X}}{\|\tilde{X}\|}, ~ Y_{n} = \varphi\left(\psi^{-1}(\|\tilde{X}_{n}\|)\right) 
\frac{\tilde{X}_{n}}{\|\tilde{X}_{n}\|}, ~n \geq 1.
\]
It is easy to see that
\[
\mathbb{P}\left(\|Y\| > a_{n} \right) =
\mathbb{P}\left( \left\|\varphi\left(\psi^{-1}(\|\tilde{X}\|)\right) 
\frac{\tilde{X}}{\|\tilde{X}\|} \right \| > \varphi(n) \right) =
\mathbb{P} \left(\|\tilde{X} \| > b_{n} \right), ~n \geq 1.
\]
It thus follows from (4.11) that $\{Y, Y_{n}; ~n \geq 1 \}$ 
is a sequence of i.i.d. symmetric {\bf B}-valued random variables such that
\[
\lim_{n \rightarrow \infty} n \mathbb{P}\left(\|Y\| > a_{n} \right) = 0
\]
and hence that, by (1.2), 
\begin{equation}
\frac{\sum_{i=1}^{n} Y_{i}}{a_{n}} \rightarrow_{\mathbb{P}} 0.
\end{equation}
By Theorem 2.1 (ii) together with (4.8) and (4.12), we have that
\[
\begin{array}{lll}
\mbox{$\displaystyle
\mathbb{P}\left(\left\|\sum_{i=1}^{n} \tilde{X}_{i} \right\| > \epsilon b_{n} \right)$}
& \leq &
\mbox{$\displaystyle  
4 \mathbb{P} \left(\left\|\sum_{i=1}^{n} \varphi\left(\psi^{-1}(\|\tilde{X}_{i}\|)\right) 
\frac{\tilde{X}_{i}}{\|\tilde{X}_{i}\|} \right\| > \epsilon a_{n} \right) 
+ \sum_{i=1}^{n}\mathbb{P}\left(\|\tilde{X}_{i}\| > b_{n} \right) $}\\
&&\\
&=& 
\mbox{$\displaystyle  
4 \mathbb{P} \left(\left\|\sum_{i=1}^{n} Y_{i} \right\| > \epsilon a_{n} \right) 
+ n \mathbb{P} \left(\|\tilde{X}\| > b_{n} \right)$} \\
&&\\
& \rightarrow &
\mbox{$\displaystyle 0 ~~\mbox{as}~ n \rightarrow \infty~~\forall~\epsilon > 0$}
\end{array}
\]
and hence that
\[
\frac{S_{n} - S_{n}^{\prime}}{2b_{n}} = \frac{\sum_{i=1}^{n} \tilde{X}_{i}}{b_{n}} 
\rightarrow_{\mathbb{P}} 0,
\]
where $S_{n}^{\prime} = \sum_{i=1}^{n} X_{i}^{\prime}$, $n \geq 1$. We thus conclude that
\[
\frac{S_{n} - S_{n}^{\prime}}{b_{n}} \rightarrow_{\mathbb{P}} 0.
\]
By Corollary 3.2, (4.8) follows. 

It remains to show that (4.9) implies (4.8). Let $\{X, X_{n}; ~n \geq 1 \}$ 
be a sequence of i.i.d. {\bf B}-valued random variables with (4.9). 
By the remarable L\'{e}vy inequality in a Banach space setting, we have that
\[
\mathbb{P} \left( \frac{\max_{1 \leq i \leq n} \left \|X_{i} - X_{i}^{\prime} \right\|}{b_{n}}
> 2 \lambda \right) 
\leq 2 \mathbb{P} \left( \frac{\left \|S_{n} - S_{n}^{\prime} \right\|}{b_{n}}
> 2 \lambda \right) 
\leq 4 \mathbb{P} \left( \frac{\left \|S_{n} \right\|}{b_{n}}
> \lambda \right) ~~\forall~n \geq 1.
\]
Then it follows from (4.9) that
\[
\lim_{n \rightarrow \infty } 
n \mathbb{P} \left( \frac{\|X - X^{\prime}\|}{b_{n}} > 2 \lambda \right)= 0, ~~\mbox{i.e.,}~
\lim_{n \rightarrow \infty } 
n \mathbb{P} \left( \left\|\frac{X - X^{\prime}}{2 \lambda} \right| > b_{n} \right)= 0.
\]
That is, (4.10) holds with $X$ replaced by symmetric random variable $(X - X^{\prime})/(2\lambda)$. 
Since (4.8) and (4.10) are equivalent, we conclude that
\[
\frac{\sum_{i=1}^{n} \frac{X_{i} - X_{i}^{\prime}}{2 \lambda}}{b_{n}} \rightarrow_{\mathbb{P}} 0,
~~\mbox{i.e.,}~ \left(\frac{1}{2 \lambda} \right) \frac{S_{n} - S_{n}^{\prime}}{b_{n}} \rightarrow_{\mathbb{P}} 0.
\]
Thus
\[
\frac{S_{n} - S_{n}^{\prime}}{b_{n}} \rightarrow_{\mathbb{P}} 0
\]
which, by Corollary 3.2, implies (4.8). ~$\Box$

\vskip 0.3cm

\noindent {\it Proof of Theorem 1.2}~~To establish this theorem, it suffices to show that, 
for every sequence $\{X, X_{n}; ~n \geq 1 \}$ of i.i.d. {\bf B}-valued random variables, the 
following three statements are equivalent:
\begin{equation}
\lim_{n \rightarrow \infty} \frac{S_{n}- \gamma_{n}}{b_{n}} 
\rightarrow 0~~\mbox{a.s.,}
\end{equation}
\begin{equation}
\limsup_{n \rightarrow \infty} \frac{\left\|S_{n}- \gamma_{n} \right\|}{b_{n}} 
< \infty~~\mbox{a.s.,}
\end{equation}
\begin{equation}
\sum_{n=1}^{\infty}\mathbb{P}(\|X\| > b_{n}) < \infty.
\end{equation}

The three statements (4.13)-(4.15) are equivalent if we can show that (4.13) and (4.15) are 
equivalent and (4.13) and (4.15) are equivalent.

For establishing the implication ``(4.13) $\Rightarrow$ (4.15)", let $\{X, X_{n};~n \geq 1 \}$ be a 
sequence of i.i.d. {\bf B}-valued random variables with (4.15). It follows from (4.15) that
\[
\lim_{n \rightarrow \infty} \frac{\sum_{i=1}^{n} \left(X_{i} - X_{i}^{\prime} \right)}{b_{n}} = 0
~~\mbox{a.s.}
\]
which implies that
\begin{equation}
\frac{X_{n} - X_{n}^{\prime}}{b_{n}} \rightarrow 0 ~~\mbox{a.s.}
\end{equation}
By the Borel-Cantelli lemma, (4.16) is equivalent to 
\[
\sum_{n=1}^{\infty} \mathbb{P} \left(\|X_{n} - X_{n}^{\prime}\| > \epsilon b_{n} \right) 
< \infty ~~\forall ~\epsilon > 0,
\]
i.e., 
\begin{equation}
\sum_{n=1}^{\infty} \mathbb{P} \left(\|X - X^{\prime}\| > \epsilon b_{n} \right) 
< \infty ~~\forall ~\epsilon > 0.
\end{equation}
Note that $\left \{\|X^{\prime}\| \leq b_{n}/2, \|X\| > b_{n} \right \} 
\subseteq \left \{ \|X - X^{\prime} \| > b_{n}/2 \right \}$ and
\[
\lim_{n \rightarrow \infty} \mathbb{P} \left(\|X^{\prime}\| \leq b_{n}/2 \right) = 1.
\]
We thus have that, for all large $n$,
\[
\mathbb{P} \left(\|X\| > b_{n} \right) \leq 2 \mathbb{P} \left \{ \|X - X^{\prime} \| > b_{n}/2 \right \} 
\]
which, together with (4.17), implies (4.15). 

We now prove ``(4.15) $\Rightarrow$ (4.13)". Let $\{X, X_{n}; ~n \geq 1 \}$ be a sequence of i.i.d. 
{\bf B}-valued random variables with (4.15). Since $\{a_{n}; n \geq 1\}$ and $\{b_{n}; n \geq 1\}$ 
are increasing sequences of positive real numbers with (1.1), by Lemma 4.1, there exist
two continuous and increasing functions $\varphi(\cdot)$ and $\psi(\cdot)$
defined on $[0, \infty)$ such that both (2.1) and (4.1) hold. Write
\[
\tilde{X} = \frac{X - X^{\prime}}{2}, ~\tilde{X}_{n} = \frac{X - X^{\prime}}{2}, ~n \geq 1
\]
and
\[
Y = \varphi\left(\psi^{-1}(\|\tilde{X}\|)\right) 
\frac{\tilde{X}}{\|\tilde{X}\|}, ~ Y_{n} = \varphi\left(\psi^{-1}(\|\tilde{X}_{n}\|)\right) 
\frac{\tilde{X}_{n}}{\|\tilde{X}_{n}\|}, ~n \geq 1.
\]
Then $\{\tilde{X}, \tilde{X}_{n}; ~n \geq 1 \}$ is
a sequence of i.i.d. symmetric {\bf B}-valued random variables such that
\begin{equation}
\sum_{n=1}^{\infty} \mathbb{P} \left(\|\tilde{X}\| > b_{n} \right) \leq
2 \sum_{n=1}^{\infty} \mathbb{P} \left(\|X\| > b_{n} \right) < \infty
\end{equation}
and $\{Y, Y_{n};~ n\geq 1\}$ is a sequence of of i.i.d. symmetric {\bf B}-valued 
random variables such that
\begin{equation}
\sum_{n=1}^{\infty} \mathbb{P} \left(\|Y\| > a_{n} \right)
= \sum_{n=1}^{\infty} \mathbb{P} \left(\| \tilde{X} \| > b_{n} \right)
< \infty.
\end{equation}
Hence, by Lemma 4.2 (i), we conclude from (4.19) and (1.5) that 
\begin{equation}
\sum_{n=1}^{\infty} \mathbb{P}\left(\left\|\sum_{i = k_{n-1}+1}^{k_{n}}
Y_{i} \right\| > \epsilon a_{k_{n}} \right)
< \infty ~~\forall~ \epsilon > 0.
\end{equation}
By Theorem 2.1 (ii), we have that, for every $n \geq 1$,
\[
\begin{array}{lll}
\mbox{$\displaystyle
\mathbb{P}\left(\left\|\sum_{i=k_{n-1}+1}^{k_{n}} \tilde{X}_{i} \right\| > \epsilon b_{k_{n}} \right)$}
& \leq &
\mbox{$\displaystyle  
4 \mathbb{P} \left(\left\|\sum_{i=k_{n-1}+1}^{k_{n}} Y_{i} \right\| > \epsilon a_{n} \right) 
+ \sum_{i=k_{n-1} + 1}^{k_{n}}\mathbb{P}\left(\|\tilde{X}_{i}\| > b_{k_{n}} \right) $}\\
&&\\
&\leq& 
\mbox{$\displaystyle  
4 \mathbb{P} \left(\left\|\sum_{i=k_{n-1}+1}^{k_{n}} Y_{i} \right\| > \epsilon a_{n} \right) 
+ \sum_{i=k_{n-1} + 1}^{k_{n}}\mathbb{P}\left(\|\tilde{X}\| > b_{i} \right)~~\forall~\epsilon > 0.$}
\end{array}
\]
It thus from (4.18) and (4.20) that
\begin{equation}
\begin{array}{ll}
& \mbox{$\displaystyle
\sum_{n=1}^{\infty} \mathbb{P}\left(\left\|\sum_{i=k_{n-1}+1}^{k_{n}} \tilde{X}_{i} \right\| > \epsilon b_{k_{n}} \right)$}\\
&\\
& \mbox{$\displaystyle  
\leq 4 \sum_{n=1}^{\infty} \mathbb{P} \left(\left\|\sum_{i=k_{n-1}+1}^{k_{n}} Y_{i} \right\| > \epsilon a_{n} \right) 
+ \sum_{n=1}^{\infty} \sum_{i=k_{n-1} + 1}^{k_{n}}\mathbb{P}\left(\|\tilde{X}_{i}\| > b_{i} \right) $}\\
&\\
& \mbox{$\displaystyle  
= 4 \sum_{n=1}^{\infty} \mathbb{P} \left(\left\|\sum_{i=k_{n-1}+1}^{k_{n}} Y_{i} \right\| > \epsilon a_{n} \right) 
+ \sum_{n=1}^{\infty} \mathbb{P}\left(\|\tilde{X}_{i}\| > b_{n} \right) $}\\
&\\
& \mbox{$\displaystyle  
< \infty~~\forall ~\epsilon > 0.$}
\end{array}
\end{equation}
By Lemma 4.2 (ii), (4.21) is equivalent to
\[
\frac{S_{n} - S_{n}^{\prime}}{2b_{n}} = \frac{\sum_{i=1}^{n} \tilde{X}_{i}}{b_{n}} \rightarrow 0~~\mbox{a.s.}
\]
Hence 
\[
\frac{S_{n} - S_{n}^{\prime}}{b_{n}} \rightarrow 0~~\mbox{a.s.}
\]
By Corollary 3.2, (4.13) follows. 

The implication ``(4.13) $\Rightarrow$ (4.14)" is obvious.

We now establish the implication ``(4.14) $\Rightarrow$ (4.13)". Let $\{X, X_{n}; ~n \geq 1 \}$ be a sequence of i.i.d. 
{\bf B}-valued random variables with (1.8). It follows from (4.14) that
\[
\limsup_{n \rightarrow \infty} \frac{\left\|\sum_{i=1}^{n} \left(X_{i} - X_{i}^{\prime} \right) \right\|}{b_{n}} 
< \infty
~~\mbox{a.s.}
\]
which implies that
\begin{equation}
\limsup_{n \rightarrow \infty} \frac{\left\|X_{n} - X_{n}^{\prime} \right\|}{b_{n}} < \infty ~~\mbox{a.s.}
\end{equation}
By the Borel-Cantelli lemma, (4.22) is equivalent to, for some constant $0 < \lambda < \infty$,
\[
\sum_{n=1}^{\infty} \mathbb{P} \left(\|X_{n} - X_{n}^{\prime}\| > \lambda b_{n} \right) 
< \infty, ~~\mbox{i.e.,}~~
\sum_{n=1}^{\infty} \mathbb{P} \left(\left\|\frac{X - X^{\prime}}{\lambda} \right \| > b_{n} \right) 
< \infty.
\]
That is, (4.15) holds with $X$ replaced by symmetric random variable $(X - X^{\prime})/\lambda$. Since (4.13) and (4.15)
are equivalent, we conclude that
\[
\left(\frac{1}{\lambda} \right) \frac{S_{n} - S_{n}^{\prime}}{b_{n}} 
= \frac{\sum_{i=1}^{n} \frac{X_{i} - X_{i}^{\prime}}{\lambda}}{b_{n}} \rightarrow 0 ~~\mbox{a.s.}
\]
Thus
\[
\frac{S_{n} - S_{n}^{\prime}}{b_{n}} \rightarrow 0 ~~\mbox{a.s.}
\]
which, by Corollary 3.2, implies (4.13). The proof of Theorem 1.2 is therefore complete. ~$\Box$

\vskip 0.2cm

\vskip 0.5cm

\noindent
{\bf Acknowledgments}\\

\noindent The research of Deli Li was partially supported
by a grant from the Natural Sciences and Engineering Research Council of
Canada and the research of Han-Ying Liang was partially supported 
by the National Natural Science Foundation of China (grant \#: 11271286).

\vskip 0.5cm

{\bf References}

\begin{enumerate}

\item Azlarov, T. A., Volodin, N. A.: Laws of large numbers for
identically distributed Banach-space valued random variables. Teor.
Veroyatnost. i Primenen. {\bf 26}, ~584-590 (1981), in Russian.
English translation in Theory Probab. Appl. {\bf 26}, 573-580
(1981).

\item de Acosta, A.: Inequalities for {\it B}-valued random vectors
with applications to the law of large numbers. Ann. Probab. {\bf 9},
157-161 (1981).

\item Feller, W.: A Limit Theoerm for Random Variables with Infinite
moments. Amer. J. Math. {\bf 68}, 257-262 (1946).
        
\item Hoffmann-J{\o}rgensen, J.: Sums of independent Banach space valued
random variables. Aarhus Univ. Preprint Series 1972/73, no. 15 (1973).

\item Hoffmann-J{\o}rgensen, J.: Sums of independent Banach
space valued random variables. Studia Math. {\bf 52}, 159-186
(1974).

\item Hoffmann-J{\o}rgensen, J.: Probability in Banach spaces. Ecole
d'Et\'{e} de Probabilit\'{e}s de St-Flour 1976. Lecture Notes in Mathematics,
vol. 598. Springer, Berlin Heidelberg 1976, pp. 1-186.

\item Kahane, J.-P.: Some random series of functions. Heath Math. Monographs,
1968. Cambridge Univ. Press, 1985, 2nd ed. 

\item Ledoux, M., Talagrand, M.:  Probability in Banach Spaces:
Isoperimetry and Processes. Springer-Verlag, Berlin (1991).

\item Li, D. L.:  A remark on the symmetrization principle for
sequences of {\it B}-valued random elements. Acta Sci. Natur.
Univ. Jilin. {\bf 1988}, 17-20 (1988), in Chinese.

\item Maurey, B., Pisier, G.: S\'{e}ries de variables al\'{e}atoires
vectorielles ind\'{e}pendantes et propri\'{e}t\'{e}s
g\'{e}om\'{e}triques des espaces de Banach. Studia Math. {\bf 58},
45-90 (1976).

\item Marcus, M. B., Woyczy\'{n}ski, W. A.: Stable measures and
central limit theorems in spaces of stable type. Trans. Amer. Math.
Soc. {\bf 251}, 71-102 (1979).

\item Mourier, E.: El\'{e}ments al\'{e}atoires dans un espace de
Banach. Ann. Inst. H. Poincar\'{e} {\bf 13}, 161-244 (1953).

\item Petrov, V. V.: Limit Theorems of Probability
Theory: Sequences of Independent Random Variables. Clarendon Press,
Oxford (1995).

\item Pisier, G.: Probabilistic methods in the geometry of Banach
spaces, in Probability and Analysis, Lectures given at the 1st 1985
Session of the Centro Internazionale Matematico Estivo (C.I.M.E.),
Lecture Notes in Mathematics, Vol. {\bf 1206}, 167-241,
Springer-Verlag, Berlin (1986).

\item Rosi\'{n}ski, J.: Remarks on Banach spaces of stable type.
Probab. Math. Statist. {\bf 1}, 67-71 (1980).

\item Woyczy\'{n}ski, W. A.: Geometry and martingales in Banach
spaces-Part II: Independent increments, in Probability on Banach
Spaces (Edited by J. Kuelbs), Advances in Probability and Related
Topics Vol. {\bf 4} (Edited by P. Ney), 267-517, Marcel Dekker, New
York (1978).

\end{enumerate}

\end{document}